\documentclass[12pt]{article}
\usepackage{amssymb}

\newcommand{\be}{\begin{equation}}
\newcommand{\ee}{\end{equation}}
\newcommand{\bea}{\begin{eqnarray}}
\newcommand{\eea}{\end{eqnarray}}
\newcommand{\barray}{\begin{array}}
\newcommand{\earray}{\end{array}}
\newcommand{\pa}{\partial}
\newcommand{\nn}{\nonumber}
\newcommand{\bitem}{\begin{itemize}}
\newcommand{\eitem}{\end{itemize}}
\newtheorem{teo}{Theorem}[section]
\newcommand{\bt}{\begin{teo}}
\newcommand{\et}{\end{teo}}
\newtheorem{Def}{Definition}[section]
\newcommand{\bd}{\begin{Def}}
\newcommand{\ed}{\end{Def}}
\newtheorem{lem}{Lemma}[section]
\newcommand{\bl}{\begin{lem}}
\newcommand{\el}{\end{lem}}
\newtheorem{prop}{Proposition}[section]
\newcommand{\bp}{\begin{prop}}
\newcommand{\ep}{\end{prop}}
\newtheorem{cor}{Corollary}[section]
\newcommand{\bc}{\begin{cor}}
\newcommand{\ec}{\end{cor}}
\newtheorem{ex}{Example}[section]
\newcommand{\bex}{\begin{ex}}
\newcommand{\eex}{\end{ex}}
\newtheorem{rem}{Remark}[section]
\newcommand{\br}{\begin{rem}}
\newcommand{\er}{\end{rem}}

\begin{document}

\begin{center}
{\Large \textbf{Deformations of Poisson structures \\ by closed
3-forms\footnote{The work was supported by the Max-Planck-Institut
f\"{u}r Mathematik (Bonn, Germany), by the Russian Foundation for
Basic Research (project no.~08-01-00464) and by the programme
``Leading Scientific Schools'' (project no. NSh-1824.2008.1).}}}
\end{center}

\medskip

\begin{center}
{\large \bf {O. I. Mokhov}}
\end{center}

\bigskip

\begin{center}
\bf {Abstract}
\end{center}

We prove that an arbitrary Poisson structure $\omega^{ij} (u)$ and
an arbitrary closed 3-form $T_{ijk} (u)$ generate the local Poisson
structure $A^{ij} (u, u_x) = M^i_s (u, u_x) \omega^{sj} (u),$ where
$M^i_s (u, u_x) (\delta^s_j + \omega^{sp} (u) T_{pjk} (u) u^k_x) =
\delta^i_j$, on the corresponding loop space. We obtain also a
special graded $\varepsilon$-deformation of an arbitrary Poisson
structure $\omega^{ij} (u)$ by means of an arbitrary closed 3-form
$T_{ijk} (u)$.

\bigskip

\begin{center}
{\large \bf {}}
\end{center}

\medskip

In this paper we prove that an arbitrary Poisson structure
$\omega^{ij} (u)$ and an arbitrary closed 3-form $T_{ijk} (u)$
generate the local Poisson structure \be A^{ij} (u, u_x) = B^i_s (u,
u_x) \omega^{sj} (u), \label{0a}\ee where \be B^i_s (u, u_x) M^s_j
(u, u_x)  = \delta^i_j,\ \ \ \ M^s_j (u, u_x) = \delta^s_j +
\omega^{sp} (u) T_{pjk} (u) u^k_x,\label{0b}\ee i.e., the matrix
operator $A^{ij} (u, u_x)$ gives the Poisson bracket \be \{I, J\} =
\int {\delta I \over \delta u^i (x)} A^{ij} (u, u_x) {\delta J \over
\delta u^j (x)} dx \label{0c}\ee on the space of functionals on the
corresponding loop space.

Let $M^N$ be an arbitrary smooth $N$-dimensional manifold with the
local coordinates $u = (u^1, \ldots, u^N)$. By the {\it loop space}
$\Omega M$ of the manifold $M^N$ we mean, in this paper, the space
of all smooth parametrized mappings of the circle $S^1$ into $M^N$,
$\gamma : S^1 \rightarrow M^N,$ $\gamma (x) = \{u^i (x)\},$ $x \in
S^1$. The tangent space $T_{\gamma} \Omega M$ of the loop space
$\Omega M$ at the point $\gamma$ consist of all smooth vector fields
$\xi = \{\xi^i, \ 1 \leq i \leq N\}$, defined along the loop
$\gamma$ with $\xi (\gamma (x)) \in T_{\gamma (x)} M$, $\forall x
\in S^1$, where $T_{\gamma (x)} M$ is a tangent space of the
manifold $M$ at the point $\gamma (x)$. All closed 2-forms
(presymplectic structures) on the loop space $\Omega M$ that are
given by matrix operators of the form $\omega_{ij} (u, u_x, \ldots,
u_{(k)})$, i.e., all closed 2-forms of the form \be \omega (\xi,
\eta) = \int_{S^1} \xi^i \omega_{ij} (u, u_x, \ldots, u_{(k)})
\eta^j dx, \label{1}\ee where $\xi, \eta \in T_{\gamma} \Omega M$,
were completely described in [1] (see also descriptions of various
differential-geometric classes of symplectic (presymplectic) and
Poisson structures in [2]--[9]).

{\bf Theorem 1 [1].} {\it A bilinear form {\rm(\ref{1})} is a closed
skew-symmetric 2-form {\rm(}a presymplectic structure{\rm)} on the
loop space $\Omega M$ if and only if \be \omega_{ij} (u, u_x,
\ldots, u_{(k)}) = T_{ijk} (u) u^k_x + \Omega_{ij} (u), \label{2}\ee
where $T_{ijk} (u)$ is an arbitrary closed 3-form on the manifold
$M^N$ and $\Omega_{ij} (u)$ is an arbitrary closed 2-form on $M^N$.}

If the matrix $\omega_{ij} (u, u_x, \ldots, u_{(k)})$ is
nondegenerate, $\det (\omega_{ij} (u, u_x, \ldots, u_{(k)})) \neq
0$, then the corresponding presymplectic form (\ref{1}), (\ref{2})
is symplectic and the inverse matrix $\omega^{ij} (u, u_x, \ldots,
u_{(k)})$, $\omega^{is} (u, u_x, \ldots, u_{(k)}) \omega_{sj} (u,
u_x, \ldots, u_{(k)}) = \delta^i_j$, gives the Poisson structure \be
\{I, J\} = \int {\delta I \over \delta u^i (x)} \omega^{ij} (u, u_x,
\ldots, u_{(k)}) {\delta J \over \delta u^j (x)} dx \label{3}\ee on
the loop space $\Omega M$, i.e., the bracket (\ref{3}) is
skew-symmetric and satisfy the Jacobi identity. Therefore Theorem 1
gives the complete description of all nondegenerate Poisson
structures on the loop space $\Omega M$ that are given by matrix
operators of the form $\omega^{ij} (u, u_x, \ldots, u_{(k)})$, i.e.,
all the nondegenerate Poisson brackets of the form (\ref{3}), $\det
(\omega^{ij} (u, u_x, \ldots, u_{(k)})) \neq 0$ (such nondegenerate
Poisson structures were studied by Astashov and Vinogradov in [9],
see also [7]--[8] and [1]--[6]). We note that if the closed 2-form
$\Omega_{ij} (u)$ is nondegenerate, $\det (\Omega_{ij} (u)) \neq 0$,
i.e, the form $\Omega_{ij} (u)$ is symplectic on $M^N$, then the
2-form (\ref{2}) is a nondegenerate form on $\Omega M$ for any
closed 3-form $T_{ijk} (u)$ on the manifold $M^N$ since it is
obvious that in this case $\det (T_{ijk} (u) u^k_x + \Omega_{ij}
(u)) \neq 0$. Thus we can define, on the loop space of an arbitrary
symplectic manifold $M^N$, the Poisson bracket \be \{I, J\} = \int
{\delta I \over \delta u^i (x)} \omega^{ij} (u, u_x) {\delta J \over
\delta u^j (x)} dx, \label{4}\ee where \be \omega^{li} (u, u_x)
(T_{ijk} (u) u^k_x + \Omega_{ij} (u)) = \delta^l_j, \label{5}\ee $I$
and $J$ being arbitrary functionals on $\Omega M$. The Poisson
bracket (\ref{4}), (\ref{5}) is a partial case of the bracket
(\ref{0a})--(\ref{0c}), namely, the case when the Poisson structure
$\omega^{ij} (u)$ is nondegenerate, $\det (\omega^{ij} (u)) \neq 0$,
$\omega^{ij} (u) = \Omega^{ij} (u),$ $\Omega^{is} (u) \Omega_{sj}
(u) = \delta^i_j,$ since \be \omega^{ij} (u, u_x) = C^i_s (u, u_x)
\Omega^{sj} (u), \ee where \be C^l_s (u, u_x) (\delta^s_j +
\Omega^{si} (u) T_{ijk} (u) u^k_x) = \delta^l_j.\label{0ba}\ee The
case of degenerate Poisson structures $\omega^{ij} (u)$, $\det
(\omega^{ij} (u)) = 0$, is much more complicated. We note that in
contrast to the case of all closed 2-forms (presymplectic
structures) of the form (\ref{1}) (Theorem 1) the problem of
description of all degenerate Poisson structures of the form
(\ref{3}) is a very complicated and unsolved problem.

{\bf Theorem 2.} {\it An arbitrary Poisson structure $\omega^{ij}
(u)$ and an arbitrary closed 3-form $T_{ijk} (u)$ give the local
Poisson bracket {\rm(\ref{0a})--(\ref{0c})}.}

First of all, we note that obviously the matrix operator $A^{ij} (u,
u_x)$ (\ref{0a}), (\ref{0b}) is skew-symmetric.

{\bf Lemma.} {\it A skew-symmetric matrix operator $A^{ij} (u, u_x)$
{\rm (\ref{0a})} gives a Poisson bracket {\rm (\ref{0c})} if and
only if the following relations hold: \be \omega^{ij} (u)
\omega^{rp} (u) {\pa M^s_r \over \pa u^i_x} = \omega^{is} (u)
\omega^{rj} (u) {\pa M^p_r \over \pa u^i_x}, \label{r1}\ee } \bea &&
\omega^{ij} (u) \omega^{rp} (u) {\pa M^s_r \over \pa u^i} -
\omega^{ij} (u) {d \over dx} \left ( {\pa M^s_r \over \pa u^i_x}
\omega^{rp} (u) \right ) + {\pa \omega^{ij} \over \pa u^r}
\omega^{rp} (u) M^s_i (u) + \nn\\ && + \omega^{is} (u) \omega^{rj}
(u) {\pa M^p_r \over \pa u^i} + {d \over dx} \left (\omega^{is} (u)
\right ) {\pa M^p_r \over \pa u^i_x} \omega^{rj} (u) + {\pa
\omega^{is} \over \pa u^r} \omega^{rj} (u) M^p_i (u) + \nn\\ && +
\omega^{ip} (u) \omega^{rs} (u) {\pa M^j_r \over \pa u^i} + {d \over
dx} \left (\omega^{ip} (u) \right ) {\pa M^j_r \over \pa u^i_x}
\omega^{rs} (u) + {\pa \omega^{ip} \over \pa u^r} \omega^{rs} (u)
M^j_i (u) = 0. \label{r2}\eea

If $M^i_s (u, u_x) = \delta^s_j + \omega^{sp} (u) T_{pjk} (u)
u^k_x$, then relations (\ref{r1}), (\ref{r2}) hold for an arbitrary
Poisson structure $\omega^{ij} (u)$ and an arbitrary closed 3-form
$T_{ijk} (u)$.

Let us add an arbitrary parameter $\varepsilon$ in the formula for
our Poisson structure:

\be A^{ij} (\varepsilon, u, u_x) = B^i_s (\varepsilon, u, u_x)
\omega^{sj} (u), \label{0aa}\ee where \be B^i_s (\varepsilon, u,
u_x) (\delta^s_j + \varepsilon \omega^{sp} (u) T_{pjk} (u) u^k_x) =
\delta^i_j.\label{0bb}\ee

We can now expand the Poisson structure $A^{ij} (\varepsilon, u,
u_x)$ in series in $\varepsilon$: \be A^{ij} (\varepsilon, u, u_x) =
\omega^{ij} (u) - \varepsilon \omega^{is} (u) T_{srk} (u)
\omega^{rj} (u) u^k_x + \cdots.\ee This expansion give an
$\varepsilon$-deformation of an arbitrary Poisson structure
$\omega^{ij} (u)$ by means of an arbitrary closed 3-form $T_{ijk}
(u)$.

We note that this $\varepsilon$-deformation of an arbitrary Poisson
structure $\omega^{ij} (u)$ belongs to a special class of graded
$\varepsilon$-deformations of Poisson structures (see, for example,
[10], [11]).

\medskip

\medskip

\medskip

{\bf {Acknowledgements.}} The work was supported by the
Max-Planck-Institut f\"ur Mathematik (Bonn, Germany), by the Russian
Foundation for Basic Research (project no.~08-01-00464) and by a
grant of the President of the Russian Federation (project no.
NSh-1824.2008.1).

\bigskip

\begin{center}
\bf {References}
\end{center}

\medskip

[1] O. I. Mokhov, ``Symplectic and Poisson geometry on loop spaces
of manifolds and nonlinear equations'', In: Topics in Topology and
Mathematical Physics, Ed. S.P.Novikov, Amer. Math. Soc., Providence,
RI, 1995, pp. 121--151;

\noindent  http://arXiv.org/hep-th/9503076 (1995).

[2] O. I. Mokhov, ``Poisson and symplectic geometry on loop spaces
of smooth manifolds'', In: Geometry from the Pacific Rim,
Proceedings of the Pacific Rim Geometry Conference held at National
University of Singapore, Republic of Singapore, December 12--17,
1994, Eds. A.J.Berrick, B.Loo, H.-Y.Wang, Walter de Gruyter, Berlin,
1997, pp. 285--309.

[3] O. I. Mokhov, ``Differential geometry of symplectic and Poisson
structures on loop spaces of smooth manifolds, and integrable
systems'', Trudy Matem. Inst. Akad. Nauk, Vol. 217, Moscow, Nauka,
1997, pp. 100--134; English translation in Proceedings of the
Steklov Institute of Mathematics (Moscow), Vol. 217, 1997, pp.
91--125.

[4] O. I. Mokhov, ``Symplectic and Poisson structures on loop spaces
of smooth manifolds, and integrable systems'', Uspekhi
Matematicheskikh Nauk, Vol. 53, No. 3, 1998, pp. 85--192; English
translation in Russian Mathematical Surveys, Vol. 53, No. 3, 1998,
pp. 515--622.

[5] O. I. Mokhov, ``Symplectic and Poisson geometry on loop spaces
of smooth manifolds and integrable equations'', Moscow--Izhevsk,
Institute of Computer Studies, 2004 (In Russian); English version:
Reviews in Mathematics and Mathematical Physics, Vol. 11, Part 2,
Harwood Academic Publishers, 2001; Second Edition: Reviews in
Mathematics and Mathematical Physics, Vol. 13, Part 2, Cambridge
Scientific Publishers, 2009.

[6] O. I. Mokhov, ``Symplectic forms on loop space and Riemannian
geometry'', Funkts. Analiz i Ego Prilozh., Vol. 24, No. 3, 1990, pp.
86--87; English translation in Functional Analysis and its
Applications, Vol. 24, No. 3, 1990, pp. 247--249.

[7] A. M. Astashov, ``Normal forms of Hamiltonian operators in field
theory'', Dokl. Akad. Nauk SSSR, Vol. 270, No. 5, 1983, pp.
363--368; English translation in Soviet Math. Dokl., Vol. 27, No. 3,
1983, pp. 685--689.

[8] A. M. Vinogradov, ``Hamiltonian structures in field theory'',
Dokl. Akad. Nauk SSSR, Vol. 241, No. 1, 1978, pp. 18--21; English
translation in Soviet Math. Dokl., Vol. 19, No. 4, 1978, pp.
790--794.

[9] A. M. Astashov and A. M. Vinogradov, ``On the structure of
Hamiltonian operators in field theory'', J. Geom. Phys., Vol. 3, No.
2, 1986, pp. 781--785; English translation in Soviet Math. Dokl.,
Vol. 27, 1983, pp. 263--287.

[10] B. Dubrovin and Y. Zhang, ``Normal forms of hierarchies of
integrable PDEs, Frobenius manifolds and Gromov - Witten
invariants''; arXiv:math/0108160 (2001).

[11] B. Dubrovin, S.-Q. Liu and Y. Zhang, ``On Hamiltonian
perturbations of hyperbolic systems of conservation laws, I:
Quasi-triviality of bi- Hamiltonian perturbations,'' Comm. Pure
Appl. Math., Vol. 59, 2006, 559 - 615; arXiv:math/0410027 (2004).

\bigskip

\begin{flushleft}
{\bf O. I. Mokhov}\\
Centre for Nonlinear Studies,\\
L.D.Landau Institute for Theoretical Physics,\\
Russian Academy of Sciences,\\
Kosygina str., 2,\\
Moscow, 117940, Russia;\\
Department of Geometry and Topology,\\
Faculty of Mechanics and Mathematics,\\
M.V.Lomonosov Moscow State University,\\
Moscow, 119992, Russia\\
{\it E-mail\,}: mokhov@mi.ras.ru; mokhov@landau.ac.ru; mokhov@bk.ru\\
\end{flushleft}

\end{document}